\documentclass[12pt]{article}
\usepackage[margin=1in]{geometry}
\usepackage{float}
\usepackage{amsmath,comment}
\usepackage{amssymb}
\usepackage{amsthm}
\usepackage{enumerate}
\usepackage{url}
\usepackage{graphicx}
\usepackage[breaklinks=true]{hyperref}
\usepackage[utf8]{inputenc}
\usepackage[T1]{fontenc}
\usepackage{xcolor, soul}

\usepackage{subcaption}
\title{On the History of the Square--and--Multiply Algorithm}

\author{
 Nuh Aydin$^{1,}$\footnote{Corresponding author: aydinn@kenyon.edu} \quad Mohammad K. Azarian$^{2}$ \quad  Omid Khormali$^{2}$ \\
 Ghaya Mtimet$^{2}$ \\[1em]
{\small $^{1}$Department of Mathematics and Statistics, Kenyon College, Gambier, OH 43022, USA}\\
{\small Email: aydinn@kenyon.edu}\\
{\small $^{2}$Department of Mathematics, University of Evansville, Evansville, IN 47722, USA}\\
{\small Emails: ma3@evansville.edu, ok16@evansville.edu, gm163@evansville.edu}
}

\date{}

\begin{document}

\maketitle

\begin{abstract}

The square-and-multiply algorithm, also known as binary 
exponentiation or repeated squaring, is a standard method 
for fast exponentiation in modern computation. Its historical 
origins, however, remain uncertain. This paper examines the 
emergence and progressive formalization of the method through 
selected primary sources. Particular attention is given to 
Jamsh\={\i}d al-K\={a}sh\={\i}'s fifteenth-century 
\textit{Mift\={a}h al-His\={a}b}, where the procedure is 
presented explicitly as a general computational method and 
claimed by al-K\={a}sh\={\i} as his own innovation. Earlier 
instances of successive squaring are identified in the works 
of al-Uql\={\i}d\={\i}s\={\i} and al-B\={\i}r\={u}n\={\i}, 
although in these cases the technique appears in particular 
calculations rather than as a fully articulated general rule. 
The earliest known antecedent is found in Ping\={a}la's 
prosodic studies in ancient India (c.\ 200 BCE), which seem 
to presuppose the conceptual basis of the method in their use 
of binary representation. As part of the historical 
development, Legendre's 1798 worked example is one of the 
earliest documented European use of the algorithm which 
appears as a subordinate step within a specific 
number-theoretic computation. The evidence suggests not a 
single continuous line of transmission, but the repeated 
independent reappearance of related procedures in distinct 
contexts. By the twentieth century, square-and-multiply became 
a special case within the broader theory of addition chains. 
By exploring this intellectual progression, this paper sheds 
some light on the historical background of an algorithm that 
is prominent in modern computation.
	\end{abstract}

\textbf{Keywords:} Square-and-multiply algorithm, Binary exponentiation, Cryptography,  Mift\={a}\d{h} al-\d{H}is\={a}b, Medieval Islamic mathematics,  Ancient Indian mathematics, Al-K\={a}sh\={\i},  Al-Uql\={\i}d\={\i}s\={\i}, Al-B\={\i}r\={u}n\={\i}, Ping\={a}la, Legendre, Addition chains, Scholz–Brauer conjecture, Modular exponentiation. \\

\section{Introduction}

The square-and-multiply algorithm is a commonly used computational technique in modern cryptography  since it is a critical element for a practical implementation of several cryptosystems. These include many well-known public key cryptosystems such as RSA, Diffie-Hellman,  ElGamal, and elliptic curve protocols for classical computers. Moreover, it  plays a key role in many primality testing algorithms, as well as quantum algorithms like Shor's algorithm to factor large integers  and solve the discrete logarithm problem \cite{shor1997}. It is also an important computational aspect of security analysis of various cryptosystems, for example in the context of side-channel attacks. Due to its prominence in cryptography  and number theory (more details in Section~\ref{sec::importance}), the algorithm is discussed in almost every book on cryptography and many books on number theory. See, for example, \cite[p.~244]{ schneier1996}, \cite[p.~71]{handbook},  \cite[pp.~69-70]{elNT}, \cite[pp.~194-198]{Yan}.  The importance of this topic motivates us to study its history.\\

Despite its prominence in modern mathematics, the historical origins of this algorithm are not fully known. The question of who first discovered or used this algorithm turns out to be more complicated than one might expect. A useful starting point is Knuth~\cite{knuth1997}, one of the most authoritative references in computer science, who provides a brief historical sketch of the algorithm (\cite[p.~461]{knuth1997}). He attributes an early appearance of the method to Ping\={a}la's \textit{Chandah\'sastra} (Rule of Metrics) before 400 AD, then Knuth credits al-Uql\={\i}d\={\i}s\={\i} of Damascus (c. 952 AD) with giving a clear discussion of how to compute $2^n$ efficiently for arbitrary $n$, and Knuth finally states that the great calculator al-K\={a}sh\={i} explicitly described the general right-to-left binary method for exponentiation in 1427 AD  (\cite[p.~462]{knuth1997}). This sketch naturally raises the question of what exactly these earlier sources contain and how they relate to each other, which is the central motivation of this work.\\
	
	A number of other sources, ranging from educational websites and blog posts to some academic publications~\cite{bhavya_vedic, booksfact_pingala, grokipedia_exponentiation, kiddle_exponentiation, prepyourmath_pingala, sahu_pingala, shah_pingala,  simplewiki_exponentiation}, also attribute the algorithm to Ping\={a}la, who lived around the second century BCE. If true, this would mean the algorithm has a much older history than is generally recognized. Knuth already states that ``this method is quite ancient'' and he says it appeared in Ping\={a}la's \textit{Chandah\'sastra}. We analyze the relevant evidence and primary sources  to clarify the historical trajectory of the square-and-multiply algorithm. Since al-K\={a}sh\={i} claims to have invented the algorithm in  the fifteenth century, we  start by looking into his work.  \\

	Al-K\={a}sh\={i}  presents the method quite explicitly in his comprehensive book \textit{Mift\={a}\d{h} al-\d{H}is\={a}b}~\cite{Miftah1, Miftah2, Miftah3}, which can be translated as \textit{Calculator's Key} or  \textit{Key to Arithmetic} (we prefer the latter translation, henceforth \textit{Mift\={a}\d{h}}). He introduces it across two rules in the algebra treatise of \textit{Mift\={a}\d{h}}~\cite{Miftah3}, first in the context of computing the sum of a geometric series, and then in a rule devoted entirely to computing large powers efficiently. Notably, al-K\={a}sh\={i} frames these contributions as rules he deduced himself, suggesting he viewed the method as his own contribution. This makes \textit{Mift\={a}\d{h} al-\d{H}is\={a}b} one of the clearest and most explicit early treatments of the algorithm that we have found in the literature. More on al-K\={a}sh\={\i}'s broader mathematical contributions and the content of \textit{Mift\={a}\d{h}} can be found in~\cite{Miftah1, Miftah2, Miftah3,  Az1, Az2, Az3, Az4, Az6, Az5, Berggren2}.\\

	It turns out that two earlier Islamic scholars, al-Uql\={\i}d\={\i}s\={\i} (920-980) and al-B\={\i}r\={u}n\={\i} (973-1048), both used techniques of successive squaring several centuries before al-K\={a}sh\={i}. Al-Uql\={\i}d\={\i}s\={\i}, who was active in the tenth century~\cite{saidan_uqlidisi}, developed a systematic procedure for finding the amount in any cell of a doubling sequence using repeated squaring. Al-B\={\i}r\={u}n\={\i}, in the next century~\cite{albiruni_chronology}, applied a similar chain of squarings to solve the classical chess problem, computing $2^{64}$ through four successive squarings. However, neither of them presented the method as a general standalone method the way al-K\={a}sh\={i} did. There was no explicitly named concept of ``algorithm'' in the fifteenth century or before, but al-K\={a}sh\={i} called it a ``rule''.  \\

    We find the appearance of the method in Europe after al-K\={a}sh\={i} in the work of Legendre at the end of the 18th century. Legendre used the method once to compute large modular exponents, however, he did not present it as a general method nor did he name it. In fact, he did not even recommend this algorithm. We further trace how, in the twentieth century, the method was absorbed into the more general theory of addition chains, initiated by Scholz in 1937 and extended by Brauer in 1939, which frames square-and-multiply as the particular addition chain determined by the binary representation of the exponent. \\
    
Historical sources spanning from ancient India to nineteenth-century Europe reveal several defining features of this method's evolution. It was employed to solve specific problems across various cultures and eras—from analyzing metrical patterns in Indian prosody to primality testing in European number theory. While the underlying algorithm remained the same, its contexts of use varied significantly. None of the authors we examined described it as a general method except for al-K\={a}sh\={i} who gave a systematic description for it and categorized  it a ``rule''. He presented it a standalone procedure in his comprehensive book, rather than as subordinate step within a specific calculation.\\

	This paper investigates these historical threads. Section~\ref{sec::importance} discusses the importance and modern applications of the algorithm. Section~\ref{sec::alkashi} examines al-K\={a}sh\={i}'s treatment of the algorithm in detail. Section~\ref{sec::earlier} looks at the earlier relevant works of al-Uql\={\i}d\={\i}s\={\i} and al-B\={\i}r\={u}n\={\i}, as well as Ping\={a}la's work on this topic. Section~\ref{sec::after} highlights some episodes of related work in Europe.  Section~\ref{sec::conclusion} offers concluding remarks with a discussion of potential future work on the subject. 

\section{The Importance of the Algorithm } \label{sec::importance}
	The algorithm known variously as \emph{square-and-multiply}, \emph{binary exponentiation}, \emph{exponentiation by squaring}, \emph{repeated squaring}, and \emph{double-and-add} in the setting of elliptic curves, is a broadly applied technique in computational number theory and modern cryptography. The core idea behind all these names is actually the same. By writing the exponent in binary and performing a sequence of squarings and conditional multiplications, the number of operations needed to compute $a^{n}$ drops from $O(n)$ to $O(\log(n))$.  Computationally, this reduction is highly significant. Since this is a number-theoretical algorithm, the size of the input is  $\log(n)$, rather than $n$. The  naive algorithm computes $a^n$ by performing the sequence of computations $a^2=a\cdot a, a^3=a\cdot a^2, \dots, a^n=a\cdot a^{n-1}.$ 
Therefore, its complexity is $O(n)$ which is exponential time because $n=2^{\log(n)}$, whereas the square-and-multiply algorithm is linear in the size of the input, hence it is an efficient algorithm. To demonstrate the practical impact of the efficient algorithm, we compared the running times of the naive versus the efficient algorithm on two modern computer algebra systems Magma and Maple. Let $n=10^{12}+27, a=3567, \text{ and }e=545566011$. The Table 1 below shows the execution times of $a^e \mod n$ on Magma and Maple using the naive method and the efficient method. The difference in computational efficiency between the two algorithms becomes even more dramatic for larger values of the parameters. 

\begin{table}[htbp]
\centering
\caption{Execution Time Comparison between Naive and Efficient Methods}
\label{tab:execution_times}
\begin{tabular}{lll}
\hline
\textbf{Computer Algebra System} & \textbf{Naive Method} & \textbf{Efficient Method} \\ \hline
Magma & 593.99 s ($\approx$ 10 mins) & 0.000 s (Instantaneous) \\
Maple & 1,693.083 s ($\approx$ 28 mins) & 0.001 s (Instantaneous) \\ \hline
\end{tabular}
\end{table}

    Knuth gives a thorough treatment of this algorithm in~\cite{knuth1997}, and Schneier's \textit{Applied Cryptography}~\cite{schneier1996}  describes it as a cornerstone of practical implementation of several public key cryptosystems. The algorithm shows up in many important areas, some of which we summarize below.\\
	
	One of the most important uses of the algorithm is in computing modular exponentiation $a^{e} \mod n$, which is the core arithmetic operation in the RSA cryptosystem~\cite{rivest1978} and the Diffie--Hellman key exchange~\cite{diffie1976}. Without this technique, RSA encryption and decryption with 2048-bit keys would simply not be feasible in practice. The algorithm is also needed in the ElGamal cryptosystem and the Digital Signature Algorithm (DSA), both of which rely on computing large modular powers efficiently.\\
	
	The algorithm also plays a central role in primality testing. Both the Fermat primality test and the Miller--Rabin primality test~\cite{miller1976, rabin1980} need to compute $a^{n-1} \pmod{n}$ for large values of $n$. Without fast modular exponentiation, generating the large primes that RSA key generation requires would not be practical. It also appears as a key subroutine in the AKS deterministic primality test~\cite{agrawal2004}.\\
	
	A modified version of the algorithm is used in elliptic curve cryptography. In the setting of elliptic curves, the group operation is point addition rather than multiplication, and a scalar multiple of a point, instead of an exponent, for a large integer, needs to be computed. The analogous technique in this setting is called \emph{double-and-add}. To compute the scalar multiple $kP$ of a point $P$ for a large positive integer $k$, one doubles and conditionally adds according to the binary expansion of $k$. This is the operation at the heart of the Elliptic Curve Diffie--Hellman (ECDH) key exchange, the Elliptic Curve Digital Signature Algorithm (ECDSA), and related protocols standardized by NIST and the NSA~\cite{hankerson2004}. As Li~et~al.\ note, the double-and-add method is directly inspired by the repeated square-and-multiply algorithm and has long been considered a fundamental technique in this area~\cite{li2020}.\\
	
	The algorithm also shows up in a perhaps unexpected place, namely quantum computing. Modular exponentiation via repeated squaring is the computational bottleneck in Shor's quantum polynomial-time factoring algorithm~\cite{shor1997}, where it has to be implemented using a circuit of reversible gates. This means the algorithm matters not just for the cryptographic systems in use today but also for the quantum algorithms that could threaten those systems in the future.\\
	
	The widespread use of the square-and-multiply algorithm has also made it an important target in security research. Kocher's timing attack~\cite{kocher1996} showed that the difference in execution time between the squaring and multiplication steps can leak private key bits from RSA implementations. That finding led to a large body of work on constant-time variants, including the Montgomery ladder algorithm, and on implementations that resist power analysis attacks~\cite{coron1999}.

\section{Al-K\={a}sh\={i}'s Work } \label{sec::alkashi}

\noindent Ghiy\={a}th al-D\={\i}n Jamsh\={\i}d Mas'\={u}d al-K\={a}sh\={\i}
(1380-1429 AD, usually referred to as  al-K\={a}sh\={\i} and sometimes as al-K\={a}shan\={\i}) was one of the most renowned mathematicians and astronomers in
Persian history, and one of the most fascinating medieval Muslim
mathematicians. \ We will very briefly reintroduce his significant works in
mathematics and astronomy.\\

\noindent One of the most notable mathematical achievements of
al-K\={a}sh\={\i} is \textit{al-Ris\={a}la al-Mu\d{h}\={\i}\d{t}\={\i}yya}
(\textquotedblleft The Treatise on the Circumference\textquotedblright), which
he completed in Arabic in 1424 AD. \ In this treatise, he correctly estimated
$2\pi$ to $9$ sexagesimal digits which is equivalent to $16$ decimal places of
accuracy \cite{Az4, Az6}. \ The second important mathematical achievement of
al-K\={a}sh\={\i} is\textit{ Ris\={a}la al-Watar wa'l-Jaib }(``The Treatise on
the Chord and Sine''), which he finalized between 1424 and 1427 AD in Arabic.
\ Al-K\={a}sh\={\i} applied Ptolemy's theorem to an inscribed quadrilateral to
obtain his famous cubic equation, and then he invented an ingenious and
quickly converging iteration algorithm to ca'lculate $\sin(1^{\circ})$ to $16$
correct decimal digits (9 correct sexagesimal places) as a root of his cubic equation.  It is remarkable that al-K\={a}sh\={\i} used both geometry and algebra to approximate $\sin(1^{\circ})$ to any desired accuracy.  Not only was this one of the earliest known approximation methods in the history of mathematics, but it was also highly creative and represented a landmark achievement in medieval algebra \cite{Az2, Az3}. The third and the most well-known mathematical work of al-K\={a}sh\={\i} written in Arabic is \textit{Meft\={a}\d{h} al-\d{H}es\={a}b}~\footnote{The title appears in two transliterations in the literature: \textit{Mift\=a\d{h} al-\d{H}is\=ab}, following standard Arabic 
romanization (as in \cite{Miftah1,Miftah2,Miftah3}), and \textit{Meftah al-hesab}, used in \cite{Az5} reflecting Persian phonological conventions. Both refer to the same work.}  (``The Calculators' Key'' or ``The Key to Arithmetic''). It
contains some of al-K\={a}sh\={\i}'s original findings, as well as his improvements on earlier  works. \ Written primarily as a
textbook, \textit{Mift\={a}\d{h}} was used for more than five
centuries, not only as a textbook, but also as an encyclopedia. It served many generations of students and professionals. \ It took al-K\={a}sh\={\i} more than 7 years to bring \textit{Mift\={a}\d{h}} to fruition, on March
2, 1427 \cite{Az5}. \textit{Mift\={a}\d{h}} has been recently translated to English in full \cite{Miftah1},\cite{Miftah2}, and \cite{Miftah3}. \\

\noindent Al-K\={a}sh\={i} was also a prominent astronomer who served as the director of Samarqand observatory which was the most advanced observatory of its time.  Al-K\={a}sh\={i} completed some of his best works while at Ulugh Beg's observatory and madrasah in Samarqand (c. 1420-1429), where he was recognized as the leading mathematician and astronomer. Some mathematical discoveries of al-K\={a}sh\={\i} were in conjunction with his various works on astronomy. The best known of his works in astronomy are:\\

\noindent(i) \ \textit{Sullam al-sam\={a}' }(``The Ladder of the Sky'' or ``The Stairway of Heaven''), completed in Arabic on March 1, 1407; (ii) \textit{Mukhtasar dar 'ilm-i hay'at }(``Compendium of Science and
Astronomy''), finished in Persian in 1410-1411 AD; (iii) \textit{Ris\={a}la
dar sharh-i \={a}l\={a}t-i rasad }(``Treatise on the Explanation of
Observational Instruments''), concluded in Persian in January 1416; (iv) \textit{Nuzha al-had\={a}iq }(``The Garden Excursion''), finalized in Arabic on February 10, 1416; and  (v) \ \textit{Z\={\i}j-i Kh\={a}q\={a}n}\={\i}
(``The Kh\={a}q\={a}n\={\i} Astronomical Tables'') was
completed in Persian 1413-1414 AD (this is a revised version of
\textit{Z\={\i}j-i Ilkh\={a}n}\={\i} of Nas\={\i}r al-D\={\i}n
al-T\={u}s\={\i}) \cite{Az1}.\\

The method of square-and-multiply appears in \textit{Mift\={a}\d{h}} even though it is not a specifically named method. It is first mentioned in the context of computing the sum of a geometric series. The fifth treatise of \textit{Mift\={a}\d{h}}  is devoted to algebra \cite{Miftah3}, which is about a third of the entire book. The other two major subjects in \textit{Mift\={a}\d{h}}   are arithmetic which consists of three treatises \cite{Miftah1},  and geometry \cite{Miftah2}.  The third chapter of Treatise 5 is titled ``On Some Computational Rules that are Crucial to Find Unknowns'' and it contains fifty rules. In Rule 9,  al-K\={a}sh\={i} describes the square-and-multiply algorithm for the first time in \textit{Mift\={a}\d{h}}. However, the main goal of this rule is to give a formula for the sum of a geometric series. Later, he gives Rule 16 which is specifically devoted to the algorithm.  Here is the full statement of Rule 9 as translated in \cite{Miftah3}:\\

\begin{quote}
{\bf The Ninth Rule.} On the sum of numbers resulting from doubling one or another number. This is also from what we deduced. The way to do this is that if the last number [of this sequence] is known, we subtract one from its double. The result is the sum of these numbers. If the last number is unknown, we look at how many times we duplicated. That is the rank of a multiple. That multiple is obtained from a first ratio of two. The way to obtain this [i.e., obtain the multiple] is to see if the rank can be halved [over and over] until it becomes one. We look then, to find how many times it can be halved until it gets to one, or we know the power of two with its rank. We square two again and again as many times as to get that number [i.e. as many as that rank]. That is, we multiply two by itself then the product by itself then the second product by itself and so on, as many times as to get that number. To obtain the last number, we double it [previous number] and always subtract one from it. The result is the sum of these numbers. If we initially add one to the numbers that are multiples,  and if the sum [of the multiples] can be halved until it becomes one, then we use it to do the same computations as earlier, and the result is the sum of the numbers with one added to it. 
\end{quote}

{\bf An Analysis of al-K\={a}sh\={i}'s Ninth Rule:} As it was custom at the time, al-K\={a}sh\={i} explains everything in prose. Algebraic symbolism was not yet introduced. It was al-Khw\={a}rizm\={i}, arguably known as the father of algebra, who wrote his algebra book in early ninth century \cite{algebra}, and that was the first time in history that algebra was specifically named as a technical mathematical operation. The English word ``algebra'' is a transliteration of the Arabic word \textit{al-jabr} that appears in the title of this book \cite{algebra}. Six centuries later, when \textit{Mift\={a}\d{h}} was written, there were still no symbols in algebra. Al-K\={a}sh\={i} explains that this rule is about finding the sum of a geometric series with the common ratio of 2 whose initial term is either 1 or some other number. In the subsequent explanations, he seems to be working with the case where the first term is 1. Therefore, he considers the geometric series $\displaystyle{\sum_{i=1}^n 2^{i-1}}$. If we know the sum of this series, then the case when the first term is some other number $a$ is easy to handle because $\displaystyle{\sum_{i=1}^n a2^{i-1}=a\sum_{i=1}^n 2^{i-1}}$. Al-K\={a}sh\={i} implies that this formula is among the things he discovered or invented. This claim needs to be examined  by researchers.  Al-K\={a}sh\={i} says if the last term of the sequence is known, then the sum equals twice the last term minus 1; that is, he gives the formula $\displaystyle{1+2+2^2+\cdots+2^n=2^{n+1}-1}$ which is the correct formula.  Then he considers the case where the last term in the series is not known. In this case, he finds the exponent (say $n$) of the last term of the series ($2^n$) from the number of terms. To compute the numerical value of the last term, he divides $n$ by $2$ consecutively until he obtains $1$, if possible. First, he assumes that $n$ can be divided by $2$ consecutively until  $1$ is obtained, in say $m$ steps. Then, he computes $2^n$ as being the $m$ successive squares  $\displaystyle \left ( \left (2^2\right )^2 \right )^{{\cdot}^{{\cdot}^{\cdot}}}$.\\

After this rule, al-K\={a}sh\={i} gives a sequence of examples. In the first example, he computes the sum $\displaystyle{\sum_{i=1}^8 2^i}$. As a first step, he computes $2^8$ via the square-and-multiply method, i.e., by computing $((2^2)^2)^2$ which requires 3 squaring operations as opposed to the direct method of computing $2,2\cdot 2, 2\cdot 4, 2\cdot 8, 2\cdot 16, 2\cdot 32, 2\cdot 64,  2\cdot 128=256$ which requires 7 multiplications.  He doubles the last number,  and subtracts 1 from it to get  $2\cdot 256-1=511$, which is $\displaystyle{\sum_{i=0}^8 2^i}$. He subtracts 1 again to find $\displaystyle{\sum_{i=1}^8 2^i}=510$.
Next, he considers the well known chessboard problem which amounts to computing $2^{63}$. The way al-K\={a}sh\={i} solves this problem is to first compute $2^{64}$ via the square-and-multiply method, that is by computing $2^2,4^2, 16^2, 256^2, 65536^2, 4294967296^2$, and $18446744073709551616^2$,  because 64 can be repeatedly halved until 1 is obtained (i.e., it is a power of 2). Then he halves the last number to obtain $2^{63}$ which he finds to be 9,223,372,036,854,775,808. After this example, al-K\={a}sh\={i} makes the following comment:

\begin{quote}
As for when the number of doubles cannot be halved [over and over] to get to one, we subtract from it the largest number that can be halved until it gets to one, then the same from the remainder, and so on until nothing is left or there is one left. So, the number is decomposed into these numbers.

\end{quote}

And he gives a quick example of this case:

\begin{quote}
{\bf Example.} If the number is ten, we decompose it into two parts which are eight and two, each of which can be halved till one.
\end{quote}

Next, he considers an example where the number is broken into three parts such that each part is a power of two. In this example, the parts are 64, 32, and 4. Therefore, the number itself is $n=100$, though al-K\={a}sh\={i} does not explicitly name 100. He gives a detailed explanation of how to compute $2^{100}$, then obtains  $\displaystyle{\sum_{i=0}^{100} 2^i}.$ In the next example, he explains how to compute $\displaystyle{\sum_{i=0}^{11} 2^i}$ in which he starts with breaking up 11 as $11=8+2+1.$ Next, al-K\={a}sh\={i} makes the point that if the first value in the sequence is a number other than 1 and if we double each term at every step, then we can find the result using the formula $\displaystyle{\sum_{i=1}^n a2^{i-1}=a\sum_{i=1}^n 2^{i-1}}$. \\

Then, in rules numbered 10 to 15,  al-K\={a}sh\={i} gives the formulas below, all in words and without proofs. He states that the rules  7, 9, 15, and 16  are his own discoveries. For applications of Rule 7 in modern mathematics, see \cite{gapSum}. Here are the formulas in modern notation that al-K\={a}sh\={i} stated in words in these six rules. We present these equations in a compact symbolic form to illustrate both the types of computations  al-K\={a}sh\={i} focused on in this section of the \textit{Mift\={a}\d{h}} and the kinds of mathematical facts he designated as `rules'.\\

\noindent\textbf{Rule 10.}%
\[
{\sum_{i=1}^{n}i(i+1)=\frac{[(n+1)-1]2}{3}\sum_{i=1}^{n+1}i=\frac
{n(n+1)(n+2)}{3}}
\]

\noindent\textbf{Rule 11.}%
\[
{\sum_{i=1}^{n}i(i+1)(i+2)=\left[  \sum_{i=1}^{n+1}i\right]  \left[  \left(
\sum_{i=1}^{n+1}i\right)  -1\right]  =\left[  \frac{(n+1)(n+2)}{2}\right]
\left[  \frac{(n+1)(n+2)}{2}-1\right]  }
\]

\noindent\textbf{Rule 12.}%
\[
{\sum_{i=1}^{n}i^{2}=\frac{2n+1}{3}\sum_{i=1}^{n}i=\frac{n(n+1)(2n+1)}{6}}
\]

\noindent\textbf{Rule 13.}%
\[
{\sum_{i=1}^{n}i^{3}=\left(  \sum_{i=1}^{n}i\right)  ^{2}=\left(
\frac{n(n+1)}{2}\right)  ^{2}}
\]

\noindent\textbf{Rule 14.}%
\[
{\sum_{i=1}^{n}i^{4}=\left[  \frac{\left(  \sum_{i=1}^{n}i\right)  -1}{5}%
+\sum_{i=1}^{n}i\right]  \sum_{i=1}^{n}i^{2}}
\]

\noindent\textbf{Remark. \ }It is intriguing that about four centuries prior
to al-K\={a}sh\={\i}, Ibn al-Haytham (c.\thinspace965 -- c.\thinspace1040),
known in the west as Alhazen, discussed 
${\displaystyle\sum_{i=1}^{n}i^{4}}$ \ in his treatise on the
\textit{Determination of the Sums of Powers of Numbers }(c. 1000--1020) \cite{Az5}. 
\bigskip

\noindent\textbf{Rule 15.}%
\[
\sum_{i=1}^{n}a^{i}=\frac{aa^{n}-a}{a-1}=\frac{a\left(  a^{n}-1\right)  }%
{a-1}=\frac{a^{n}-a}{a-1}+a^{n},
\]

\noindent where $a$ is any number except $1$. \ For $a=\displaystyle\frac
{p}{q}<1$, al-K\={a}sh\={\i} presented the following formula%
\[
\sum_{i=1}^{n}\left(\frac{p}{q}\right)^{i}=\frac{(q^{n}-p^{n})p}{(q-p)q^{n}}.
\]

Then comes the sixteenth rule which is entirely devoted to the square-and-multiply method. Here is the full statement of the rule given in \cite{Miftah3}.

\begin{quote}
\textbf{The Sixteenth Rule.} To get a large power of a number without finding the consecutive powers between them. This is also from what we deduced. We take the power and check if it can be halved until one. If it is, we find the number of times it needs to be halved to get to one. Then, we square the first power that many times to get the final power, which is the wanted value.
\end{quote}

Then al-K\={a}sh\={\i}  gives two examples. In the first example, he computes $5^8$. Here, 8 is a power of 2 and 1 is obtained by halving 8 three times. Therefore, he obtains the result by 3 squarings: $5^2=25, 25^2=625$ and $625^2=390625$.  He then explains how to apply the method if the exponent is not a power of 2. He writes:

\begin{quote}
If the power cannot be halved until one, we subtract from it the largest number possible that can be halved to one, then do the same with the remainder, and so on until either nothing is left or one is left. We get numbers whose sum equals the wanted power, and each of them can be halved until one or one of them is one itself and the others can be halved until one. We put the numbers in a table as mentioned earlier in the ninth rule. We find the count of halvings required for each number to get to one, and write it next to its respective number. We write a zero next to the one. We call these counts the number-of-times. Then, we square the first power again and again as many times as the largest of the numbers of times. We write the last square next to its respective number. Similarly, we write next to each number the last square which results from squaring the first power a number of times corresponding to that number. We put the first power next to the zero. Then, we multiply these powers in the column by each others, the result is the last power that is wanted.
\end{quote}

In this case, al-K\=ash{\=\i}  explains the technique to compute  $a^n$ where $n$ is not a power of 2 as follows. He breaks up $n$ as $n=n_1+n_2+\cdots +n_k$ where the numbers $n_1, n_2, \ldots,\, n_k$  are powers of 2, i.e., he obtains the binary (base-2) representation  of $n$. Then, he computes $b_i=a^{n_i}$ for each $i$, using the 16th rule.   He then multiplies all $b_i$'s to obtain $a^n$. Next, he gives an example  to illustrate how to apply this technique to compute $4^{14}$. First he writes $14=8+4+2$. Then he computes $4^8, 4^4,$ and $4^2$, using the rule, and finally multiplies all of the numbers together to obtain $4^{14}$. 

\section{Earlier Works  } \label{sec::earlier}
As noted in the introduction, two Islamic scholars predate al-K\={a}sh\={i} in their use of successive squaring. Al-Uql\={i}d\={i}s\={i}, in the tenth century, and al-B\={i}r\={u}n\={i}, in the eleventh century, both employed repeated squaring to compute large powers of~$2$ efficiently, several centuries before al-K\={a}sh\={i} articulated the method as a general rule. However, neither of them  presented the technique as a standalone general method. Rather, in both cases the method was embedded within specific computational problems. In this section, we examine what each of these scholars did, followed by a discussion of Ping\={a}la's  earlier work.

\subsection{Al-B\={i}r\={u}n\={i} (973--1048 CE)}

In \textit{The Chronology of Ancient Nations}, 
al-B\={i}r\={u}n\={i} addresses the same classical 
chess problem, requiring the computation of 
$2^{64} - 1$ ~\cite[~pp.~132--136]{albiruni_chronology}. He presents 
two fundamental rules.\\

The first rule states that the square of the number 
of a check $x$ of the 64 checks of the chessboard 
is equal to the number of that check whose distance 
from check $x$ is equal to the distance of check $x$ 
from the first check. As an example, al-B\={i}r\={u}n\={i} 
takes the square of the number of the 5th check, 
i.e., \ $16^2 = 256$, which is the number belonging 
to the 9th check, since the distance of the 9th check 
from the 5th equals the distance of the 5th check 
from the first ~\cite[~p.~134]{albiruni_chronology}.\\

The second rule states that the number of a check $x$ 
minus $1$ is equal to the  total sum of the numbers 
of all the preceding checks. As an example, the number 
of the 6th check is $32$, and $32 - 1 = 31$, which 
equals $1 + 2 + 4 + 8 + 16$ ~\cite[~p.~132]{albiruni_chronology}.\\

Using the first rule, al-B\={i}r\={u}n\={i} computes 
$2^{64}$ through a chain of successive squarings. 
He expresses this as~\cite[~p.~132]{albiruni_chronology}:
\[
\left\{\left[\left(16^2\right)^2\right]^2\right\}^2 
= 16^{16} = 2^{64},
\]
tracing the following chain of checks and their values:

\begin{table}[h]
	\centering
	\begin{tabular}{ccc}
		\hline
		\textbf{Check Number} & \textbf{Operation} & 
		\textbf{Value} \\
		\hline
		5th  & starting point  & $16 = 2^4$   \\
		9th  & square of 5th     & $256 = 2^8$  \\
		17th & square of 9th     & $65{,}536 = 2^{16}$ \\
		33rd & square of 17th    & $4{,}294{,}967{,}296 = 2^{32}$ \\
		65th & square of 33rd    & 
		$18{,}446{,}744{,}073{,}709{,}551{,}616 = 2^{64}$ \\
		\hline
	\end{tabular}
	\caption{Al-B\={i}r\={u}n\={i}'s successive squaring 
		chain for the chess problem~\cite[~pp.~134--135]{albiruni_chronology}.}
	\label{tab:biruni_squaring}
\end{table}

He then subtracts 1 to obtain the total sum of all 
the checks of the chessboard, arriving at 
$18{,}446{,}744{,}073{,}709{,}551{,}615$ ~\cite[~p.~132]{albiruni_chronology}.

\subsection{Al-Uql\={\i}d\={\i}s\={\i} (c.\ 920--980 CE)}

In \textit{The Arithmetic of al-Uql\={\i}d\={\i}s\={\i}}, Book~IV, Chapter~32, titled ``On Doubling One, Sixty-Four Times,'' al-Uql\={\i}d\={\i}s\={\i} addresses the problem of computing the number in any cell of a doubling sequence, as well as the total of all cells up to a given one~\cite[pp.~337--342]{saidan_uqlidisi}.

Al-Uql\={\i}d\={\i}s\={\i} begins by establishing a 
\textit{standard of ten cells}:
\[
1 \quad 2 \quad 4 \quad 8 \quad 16 \quad 32 \quad 
64 \quad 128 \quad 256 \quad 512,
\]
to which he reduces all subsequent computations ~\cite[p.~339]{saidan_uqlidisi}. He observes that multiplying the number in any cell by its like (i.e., \ squaring it) gives the number in twice the cell number less 1.  For example, the number on the fifth cell is 16. If you square 16, you get 256 which  is the number on the ninth cell, and $2\cdot5 -1=9$. The total of all cells up to the last is obtained by doubling the last number and subtracting 1, i.e.,
\[
\sum_{i=1}^{n} 2^{i-1} = 2^n - 1.
\]

Al-Uql\={\i}d\={\i}s\={\i}'s method for finding the number in 
an odd-numbered cell proceeds by adding 1 to the cell number and halving it, then using the ten-cell standard. For the 17th cell, he states: ``we add one to seventeen and take half of that, finding nine. 
We then multiply what is in the ninth cell by its like: $256 \times 256 = 65{,}536$, which is the amount in the 17th cell.''\\

For the 51st cell, al-Uql\={\i}d\={\i}s\={\i} proceeds as 
follows ~\cite[pp.~341--342]{saidan_uqlidisi}: increase 51 by 1
and take half, which is 26; take its half, which is 13; add 1 and take half, which is 7; add 1 and take half, which leads to the fourth cell, in which we find 8. He then computes upward through successive 
squarings:
\begin{align*}
	8^2 &= 64 \quad \text{(7th cell)}, \\
	64^2 &= 4{,}096 \quad \text{(13th cell)}, \\
	4{,}096^2 &= 16{,}777{,}216 \quad \text{(25th cell)}, \\
	2 \times 16{,}777{,}216 &= 33{,}554{,}432 
	\quad \text{(26th cell)}, \\
	33{,}554{,}432^2 &= 1{,}125{,}899{,}906{,}842{,}624 
	\quad \text{(51st cell)}.
\end{align*}
That is, after reaching the 25th cell, he doubles 
the amount to reach the 26th cell, and then squares 
it to reach the 51st cell, since $2 \times 26 - 1 = 51$.\\

For even-numbered cells, al-Uql\={i}d\={i}s\={i} states that we halve the cell number and add 1, then proceed as in the odd case~ \cite[p.~342]{saidan_uqlidisi}. For example, consider the 14th cell. Take half of 14, which is 7. Then, add 1 and take half of that, which leads to the fourth cell, in which we find 8. He then computes
\begin{align*}
8^2 &= 64 \quad \text{(7th cell)}, \\
64^2 &= 4{,}096 \quad \text{(13th cell)}, \\
2 \times 4{,}096 &= 8{,}192 \quad \text{(14th cell)},
\end{align*}
that is, after reaching the 13th cell by two successive squarings, he doubles the amount to obtain the value in the 14th cell.

\subsection{Ping\={a}la (c.\ 2nd century BCE)}
    Ping\={a}la is an ancient Indian scholar who is credited with writing the \textit{Chandah\'{s}utra}, a treatise on Sanskrit prosody concerned with the systematic study of metrical patterns in poetry. A number of sources, ranging from educational websites and blog posts to some academic publications~\cite{simplewiki_exponentiation, kiddle_exponentiation, prepyourmath_pingala, booksfact_pingala, bhavya_vedic, shah_pingala, sahu_pingala}, attribute the square-and-multiply algorithm to Ping\={a}la, sometimes calling it the earliest known instance of binary exponentiation. Knuth states in~\cite{knuth1997} that the binary method appeared before 400 AD in Ping\={a}la's \textit{Chandah\'{s}utra}. Therefore, Ping\={a}la's work  deserves a careful examination in any discussion of the history of the square-and multiply-algorithm.\\
	
	One thing we noticed is that many of the sources  do not directly  engage with Ping\={a}la's original text and appear to rely on each other rather than on independent primary source analysis. For instance, the blog posts and educational websites~\cite{prepyourmath_pingala, booksfact_pingala, bhavya_vedic, simplewiki_exponentiation, kiddle_exponentiation} make statements that  Ping\={a}la invented binary exponentiation without providing  direct textual evidence from the \textit{Chandah\'{s}utra} itself. The academic paper by Sahu~\cite{sahu_pingala} similarly attributes algorithmic thinking to Ping\={a}la, but draws on secondary literature rather than primary sources.\\
	
	One of the scholarly treatments of Ping\={a}la's original work that we found is the paper by Shah~\cite{shah_pingala}, which systematically traces Ping\={a}la's algorithm through Indian mathematical literature over more than a millennium. Shah does describe a divide-and-conquer recursive algorithm in the \textit{Chandah\'{s}utra} for computing powers of 2, specifically the algorithm for calculating the total number of possible metrical forms of a given length. In Ping\={a}la's scheme, if the number of syllables $n$ can be halved, the result is obtained by squaring, and if not, we subtract 1 and double. This is expressed in four cryptic s\={u}tras and is used to compute $S_n = 2^n$, the total count of all forms of an $n$-syllable meter. The algorithm is recursive and it does rest on the same underlying mathematical principle as the square-and-multiply technique. In this sense, Ping\={a}la does use a method that is mathematically equivalent to binary exponentiation.\\
	
	  It is interesting to note that Ping\={a}la's use of this technique is entirely embedded within the specific problem of counting metrical patterns in Sanskrit poetry. It is not presented as a general computational method for computing arbitrary powers of a given number. Joseph's \textit{The Crest of the Peacock}~\cite{joseph2010crest} discusses the \textit{Chandah\'{s}utra} exclusively in the context of combinatorics and Pascal's triangle, with no mention of a general exponentiation algorithm. \\
      
      A well known expert  in the history of Indian mathematics is Kim Plofker. Plofker's chapter in~\cite{plofker_katz} and her book \cite{Plofker2009} describe Ping\={a}la's work related to binary exponentiation in the context of prosody and combinatorics . In the section titled ``Mathematical Ideas in Other Disciplines'' in her book, Plofker describes Ping\={a}la's method of using a succession of squaring and doubling operations  to calculate $2^n$. She (rather than Ping\={a}la himself) illustrates the algorithm  via the example of  computing $2^7$.  In section 5.3 of \cite{Plofker2009}, Plofker states  that a later Indian scholar presents a more general version of Ping\={a}la's technique to compute $r^n$. Plofker states that this scholar was a medieval (ninth century) namesake of Mahavira, the founder of Jainism (6th-5th century BCE).
      
\section{Developments after al-K\=ash\=\i}\label{sec::after}

When we follow the algorithm forward from al-K\=ash\=\i, it appears in Europe about four centuries later. Knuth's account of the binary method \cite[pp.~461--462]{knuth1997} traces it to Ping\=ala, al-Uql\=\i d\=\i s\=\i, al-B\=\i r\=un\=\i, and al-K\=ash\=\i, and does not mention Legendre. Our investigation indicates that Adrien-Marie Legendre (1752-1833) employs the method in a worked example in his \emph{Essai sur la th\'eorie des nombres} (Paris, 1798, p.~229) \cite{Legendre1798}. As in the earlier works of al-Uql\=\i d\=\i s\=\i, al-B\=\i r\=un\=\i, and Ping\=ala, Legendre does not give the method a name. It appears instead as one step inside a specific task, namely deciding whether a prime $c$ divides $x^2 + a$. By Euler's criterion, this requires the value of the power residue $a^{(c-1)/2} \pmod{c}$. It is interesting to note that Knuth cites this very page of Legendre \cite[p.~326]{knuth1997}, but only for a different technique given there, namely converting integers to binary by repeated division, without commenting on the square-and-multiply computation on the same page.\\

The clear description of the method is in the Second Part, \S\,VII, Art.~183 (p.~229) \cite{Legendre1798}, inside the worked example that tests whether $1013$ divides $x^2 + 601$. To do the test by direct computation, Legendre has to find $601^{506} (\bmod \   1013)$, where $506 = (1013-1)/2$. He notes that $506$ in binary is $111111010$, so it is the sum of the powers of 2 with exponents $8,7,6,5,4,3,1$. Forming the powers of $601$ whose exponents are these powers of 2 takes eight squarings, and combining the seven chosen powers takes six more multiplications, for fourteen multiplications in all, with the same number of reductions modulo $1013$. He then carries out the steps in full. He computes $601^{2}, 601^{4}, 601^{8}, \dots, 601^{256}$ by repeated squaring, reducing at each step, and then he combines them as $601^{384}, 601^{448}, 601^{480}, 601^{496}, 601^{504}, 601^{506}$, ending at $-1$. This is the essentially modern square-and-multiply algorithm. It uses the binary representation of the exponent, it splits the work into a run of squarings followed by chosen multiplications, and it reduces the result modulo 1013 at each step. A footnote on the same page even gives a quick way to turn a number into binary by repeated division, shown on the example $11183445$, so Legendre treats the binary form as a standard working tool.\\

It is also worth noting the role square-and-multiply plays in Legendre's text. He does not present it as a method to adopt, but works it out in full as the laborious route that his reciprocity method avoids, showing the computation so that the reader can compare the two and see how much faster his approach is. In what may be the first clear use of binary exponentiation in Europe, the method thus appears only as a foil. While al-K\=ash\=\i \ set out square-and-multiply as a general rule and claimed it as his own discovery, Legendre never names it and uses it inside an example as the slower alternative to a method he prefers.\\

 In the twentieth century the method was absorbed into a broader theory, that of addition chains. An \emph{addition chain} for $n$ is a sequence $1 = a_0, a_1, \dots, a_r = n$ in which each term is the sum of two earlier terms. Addition chains can be used for addition-chain exponentiation so that a chain of length $r$ for $n$ yields a computation of $x^{n}$ in $r$ multiplications.  As Knuth explains, finding the most economical way to compute $x^{n}$ reduces to addition because the exponents are additive \cite[\S4.6.3, p.~465]{knuth1997}. In these terms the square-and-multiply method corresponds to the particular addition chain determined by the binary representation of the exponent. It is known that the shortest addition-chain algorithm requires no more multiplications than the square-and-multiply method and usually fewer. When Arnold Scholz introduced addition chains in 1937, in the form of a problem (Aufgabe~253) \footnote{\emph{Aufgabe} is German for ``problem''; Scholz posed his result as Problem~253 in the journal's ``Aufgaben und L\"osungen'' (Problems and Solutions) section.} \cite{Scholz1937}, he defined the minimal chain length $\ell(n)$ and stated the application directly, namely to form the $n$th power of a residue modulo $m$ with as few multiplications as possible, which is exactly the fast modular exponentiation that square-and-multiply performs. In the same problem he posed the inequality $\ell(2^{n}-1) \le n - 1 + \ell(n)$,\footnote{Scholz stated this in the equivalent form $\ell(2^{m+1}-1) \le m + l(m+1)$; the two differ only by re-indexing.} now known as the Scholz--Brauer conjecture and, to the best of our knowledge, is still open \cite[\S4.6.3, p.~475]{knuth1997}. Alfred Brauer took up the problem two years later \cite{Brauer1939}, opening his paper with the same application of computing $c^{n} \bmod m$ with the fewest multiplications. He proved that $\ell(n) \sim \log_2 n$ and introduced the restricted chains now called star chains, for which he established the corresponding inequality. Brauer's proof of the bound $l(n) \le 2m$ (for $2^m + 1 \le n \le 2^{m+1}$) is itself the square-and-multiply construction, since he writes $n$ in binary, forms $1, 2, 2^2, \dots, 2^m$ by repeated doubling, and combines the required powers, so that square-and-multiply reappears within his work as the basic chain he then sets out to improve \cite[p.~737]{Brauer1939}. This material received its standard modern treatment, together with the historical notes on which much of this section draws, in Knuth's book \cite[\S4.6.3]{knuth1997}, and from there the algorithm passed into its central role in modern cryptography and computational number theory.\\

 These later instances corroborate the trajectory observed in the earlier sources. While the underlying idea appears in the works of scholars from Ping\=ala through al-Uql\=\i d\=\i s\=\i\ and al-B\=\i r\=un\=\i, and again in Legendre, in these cases it remains a localized tool for a specific problem rather than a standalone method. Al-K\=ash\=\i's fifteenth-century treatment, by contrast, states it as a general rule. Based on the available sources, this synthesis was an important step towards a general algorithm.

\section{Conclusion and Future Work} \label{sec::conclusion}
    
The number theoretical algorithm known as  the square-and-multiply, binary exponentiation, or repeated squaring is an important algorithm in modern computational number theory and cryptography that enables fast exponentiation of numbers of the form $a^m \mod n$ for large integers $m$ and $n$. Tracing historical origins of this algorithm, we find an explicit description of the algorithm as a ``rule'' in al-K\={a}sh\={\i}'s fifteenth century book  \textit{Mift\={a}\d{h}
al-\d{H}is\={a}b}.  Although  al-K\={a}sh\={\i}'s text implies that he discovered the algorithm himself, the underlying procedure had been detailed by scholars well before  al-K\={a}sh\={\i}.  For example, two Islamic scholars   al-Uql\={i}d\={i}s\={i} and al-B\={i}r\={u}n\={i}, who predate al-K\={a}sh\={i}  by several centuries,   employed successive squaring to compute large powers  of $2$ efficiently.  Al-Uql\={i}d\={i}s\={i}'s treatment  is presented as a general procedure, while al-B\={i}r\={u}n\={i}'s  treatment is presented specifically in the context  of the chess problem, tracing a fixed chain of four  squarings from the 5th to the 65th check. Going further in history than Islamic scholars, there is evidence  that the ancient Indian scholar Ping\={a}la had the basic idea of the square-and-multiply algorithm, which is to use the binary representation of the exponent, even though it was in the context of counting metrical patterns. Using the binary representation of the exponent is the essence of the algorithm that makes it efficient compared to the naive algorithm, which has exponential running time. Although it was not initially stated as a general mathematical procedure, on the basis of available evidence, Ping\={a}la may be viewed as the first person who used a procedure that we now call the square-and-multiply algorithm. While many authors between Pingala and Legendre employed the algorithm in the context of solving specific problems, none of them, except for al-K\={a}sh\={\i}, described it as a general procedure.  \\

 This paper has primarily examined the development of the square-and-multiply algorithm within Indian and medieval Islamic mathematical traditions, with a preliminary look at later European developments. Whether binary exponentiation appears in European mathematical texts between al-K\={a}sh\={\i}'s time and Legendre's 1798 example remains to be investigated. A comprehensive investigation of European sources from this period, spanning primary mathematical texts across multiple languages and genres, is beyond the scope of this study. Furthermore, additional relevant sources may yet be uncovered within medieval Islamic and Sanskrit mathematical literature. Legendre's 1798 example shows that this gap is not entirely empty: at least one leading mathematician was using square-and-multiply as a familiar computational tool at the end of the eighteenth century. This suggests, though a single instance cannot establish it, that binary exponentiation may have been more widely known among mathematicians of Legendre's time, and possibly earlier. An interesting question for future research therefore remains: what further explicit descriptions or uses of the square-and-multiply algorithm can be identified in the centuries between al-K\={a}sh\={\i}'s \textit{Mift\={a}\d{h} al-\d{H}is\={a}b} (1427) and its eventual ubiquity in modern computational mathematics? At present, this historical gap remains a significant and open invitation for further scholarly inquiry.

\section*{Declaration of Generative AI and AI-assisted technologies in the Manuscript Preparation Process}

During the preparation of this work, the authors used a few different AI tools including chatGPT, Gemini, and Perplexity in order to find relevant publications. The authors also used these tools to do sentence level editing of certain parts of the manuscript. After using these tools, the authors reviewed and edited the content as needed and they take full responsibility for the content of the published article.

\section*{Acknowledgments}
	
	The third and fourth authors gratefully acknowledge support from the Global Scholar Award, University of Evansville Center for Innovation and Change.

\end{document}